# Multi-armed bandit problem with precedence relations

**Hock Peng Chan**[1,*], **Cheng-Der Fuh**[2,†] **and Inchi Hu**[3,‡]

*National University of Singapore, National Central University, Academia Sinica and Hong Kong University of Science and Technology*

**Abstract:** Consider a multi-phase project management problem where the decision maker needs to deal with two issues: (a) how to allocate resources to projects within each phase, and (b) when to enter the next phase, so that the total expected reward is as large as possible. We formulate the problem as a multi-armed bandit problem with precedence relations. In Chan, Fuh and Hu (2005), a class of asymptotically optimal arm-pulling strategies is constructed to minimize the shortfall from perfect information payoff. Here we further explore optimality properties of the proposed strategies. First, we show that the efficiency benchmark, which is given by the regret lower bound, reduces to those in Lai and Robbins (1985), Hu and Wei (1989), and Fuh and Hu (2000). This implies that the proposed strategy is also optimal under the settings of aforementioned papers. Secondly, we establish the super-efficiency of proposed strategies when the bad set is empty. Thirdly, we show that they are still optimal with constant switching cost between arms. In addition, we prove that the Wald's equation holds for Markov chains under Harris recurrent condition, which is an important tool in studying the efficiency of the proposed strategies.

## 1. Introduction

Suppose there are $\mathcal{U} = J_1 + \cdots + J_I$ statistical populations, $\Pi_{11}, \Pi_{12}, \ldots, \Pi_{IJ_I}$. Pulling arm $ij$ once corresponds to taking an observation from population $\Pi_{ij}$. The observations from $\Pi_{ij}$ form a Markov chain on a state space $D$ with transition probability density function $p_{ij}(x, y, \theta)$ with respect to a $\sigma$-finite measure $Q$, where $\theta$ is an unknown parameter belonging to a parameter space $\Theta$. The stationary probability distribution for the Markov chain exists and has probability density function $\pi_{ij}(\cdot, \theta)$.

At each step, we are required to sample one of the statistical populations obeying the partial order $ij \preceq i'j' \Leftrightarrow i \leq i'$. An adaptive policy is a sampling rule that dictates, at each step, which population should be sampled based on observations before that step. We can represent a policy as a sequence of random variables $\phi = \{\phi_t | \phi_{t-1} \preceq \phi_t, t = 1, 2, \ldots\}$ taking values in $\{ij | i = 1, \ldots, I; j = 1, \ldots, J_i\}$ such that the event $\{\phi_t = ij\}$ 'take an observation from $\Pi_{ij}$ at step $t$' belongs to the $\sigma$-field generated by $\phi_1, X_1, \ldots, \phi_{t-1}, X_{t-1}$, where $X_t$ denotes the state of the population being sampled at $t$-th step.

---

[1]National University of Singapore, e-mail: `stachp@nus.edu.sg`
[2]Academia Sinica, e-mail: `stcheng@stat.sinica.edu.tw`
[3]Hong Kong University of Science and Technology, e-mail: `imichu@ust.hk`
[*]Research supported by grants from the National University of Singapore.
[†]Research partially supported by the National Science Council of ROC.
[‡]Research partially supported by Hong Kong Research Grant Council.

Let the initial state of $\Pi_{ij}$ be distributed according to $\nu_{ij}(\cdot;\theta)$. Throughout this paper, we shall use the notation $E_\theta$ ($P_\theta$) to denote expectation (probability) with respect to the initial distribution $\nu_{ij}(\cdot;\theta)$; similarly, $E_{\pi(\theta)}$ to denote expectation with respect to the stationary distribution $\pi_{ij}(\cdot;\theta)$. We shall assume that $\mathcal{V}_{ij} = \{x \in D : \nu_{ij}(x;\theta) > 0\}$ does not depend on $\theta$ and $v_{ij} := \inf_{x \in \mathcal{V}_{ij}} \inf_{\theta,\theta' \in \Theta} [\nu_{ij}(x;\theta)/\nu_{ij}(x;\theta')] > 0$ for all $i,j$. Suppose that $\int_{x \in D} |g(x)| \pi_{ij}(x;\theta) Q(dx) < \infty$. Let

$$\mu_{ij}(\theta) = \int_{x \in D} g(x) \pi_{ij}(x;\theta) Q(dx)$$

be the mean reward under stationary distribution $\pi_{ij}$ when $\Pi_{ij}$ is sampled once. Let $N$ be the total sample size from all populations, and

$$(1.1) \qquad T_N(ij) = \sum_{t=1}^{N} \mathbf{1}_{\{\phi_t = ij\}}$$

be the sample size from $\Pi_{ij}$ and $\mathbf{1}$ denotes the indicator function. It follows that the total reward equals

$$(1.2) \qquad W_N(\theta) := \sum_{t=1}^{N} \sum_{i=1}^{I} \sum_{j=1}^{J_i} E_\theta\{E_\theta[X_t \mathbf{1}_{\{\phi_t = ij\}} | \mathcal{F}_{t-1}]\}.$$

In the case of independent rewards, that is, when $p_{ij}(x,y,;\theta) = p_{ij}(y;\theta)$ for all $i,j,x,y$ and $\theta$, $W_N(\theta) = \sum_{i=1}^{I} \sum_{j=1}^{J_i} \mu_{ij}(\theta) E_\theta T_N(ij)$. We shall show in the Appendix that for Markovian rewards, under regularity conditions A3–A4 (see Section 2.1), there exists a constant $C_0 < \infty$ independent of $\theta \in \Theta$, $N > 0$ and the strategy $\phi$ such that

$$(1.3) \qquad \left| W_N(\theta) - \sum_{i=1}^{I} \sum_{j=1}^{J_i} \mu_{ij}(\theta) E_\theta T_N(ij) \right| \leq C_0.$$

In light of (1.3), maximizing $W_N(\theta)$ is asymptotically equivalent [up to a $O(1)$ term] to minimizing the regret

$$(1.4) \quad R_N(\theta) := N \mu^*(\theta) - W_N(\theta) = \sum_{ij: \mu_{ij}(\theta) < \mu^*(\theta)} [\mu^*(\theta) - \mu_{ij}(\theta)] E_\theta T_N(ij),$$

where $\mu^*(\theta) := \max_{1 \leq i \leq I} \max_{1 \leq j \leq J_i} \mu_{ij}(\theta)$.

Because adaptive strategies $\phi$ that are optimal for all $\theta \in \Theta$ and large $N$ in general do not exist, we consider the class of all (asymptotically) *uniformly good* adaptive strategies under the partial order constraint $\preceq$, satisfying

$$(1.5) \qquad R_N(\theta) = o(N^\alpha), \quad \text{for all } \alpha > 0 \text{ and } \theta \in \Theta.$$

Such strategies have regret that does not increase too rapidly for any $\theta \in \Theta$. We would like to find a strategy that minimizes the increasing rate of the regret within the class of uniformly good adaptive strategies under the partial order constraint $\preceq$.

The rest of the article is organized as follows. In Section 2, we present the assumptions and introduce the concept of bad sets. The regret lower bound is investigated in Section 3. We also prove that the regret lower bound specializes to other lower bounds obtained by previous authors under less general settings. Section 4 contains the super efficiency result when the bad sets are empty. The optimality of the proposed strategies under constant switching cost is investigated in Section 5. The last section includes the proof of Wald's equation for Markov random walks under Harris recurrence condition.

## 2. The assumption and bad sets

Denote the Kullback-Leibler information number by

$$(2.1) \quad I_{ij}(\theta, \theta') = \int_{x \in D} \int_{y \in D} \log \left[ \frac{p_{ij}(x, y; \theta)}{p_{ij}(x, y; \theta')} \right] p_{ij}(x, y; \theta) \pi_{ij}(x; \theta) Q(dy) Q(dx).$$

Then, $0 \leq I_{ij}(\theta, \theta') \leq \infty$. We shall assume that $I_{ij}(\theta, \theta') < \infty$ for all $i, j$ and $\theta, \theta' \in \Theta$. Let $\mu_i(\theta) = \max_{1 \leq j \leq J_i} \mu_{ij}(\theta)$ be the largest reward in the $i$-th group of arms, and

$$(2.2) \quad \Theta_i = \{\theta \in \Theta : \mu_i(\theta) > \mu_{i'}(\theta) \text{ for all } i' < i \text{ and } \mu_i(\theta) \geq \mu_{i'}(\theta) \text{ for all } i' \geq i\}$$

be the set of parameter values such that the first optimal job is in group $i$. Let

$$(2.3) \quad \Theta_{ij} = \{\theta \in \Theta_i : \mu_{ij}(\theta) = \mu_i(\theta)\}$$

be the parameter set such that arm $ij$ is one of the first optimal ones. Each $\theta \in \Theta$ belongs to exactly one $\Theta_i$ but may belong to more than one $\Theta_{ij}$. Let

$$(2.4) \quad \Theta_i^* = \{\theta \in \Theta : \mu_i(\theta) > \mu_{i'}(\theta) \text{ for all } i' \neq i\}$$

be the parameter set in which all the optimal arms lie in group $i$. Clearly, $\Theta_i^* \subset \Theta_i$ but the reverse relation is not necessarily true.

### 2.1. The assumptions

We now state a set of assumptions that will be used to prove the optimality results. Let $\Theta$ be a compact subset of $\mathbf{R}^d$ for some $d \geq 1$.

A1. $\mu_{ij}(\cdot)$ are finite and continuous on $\Theta$ for all $i, j$. Moreover, no arm group is redundant in the sense that $\Theta_i^* \neq \emptyset$ for all $i = 1, \ldots, I$.
A2. $\sum_{j=1}^{J_1} I_{1j}(\theta, \theta') > 0$ for all $\theta' \neq \theta$ and $\inf_{\theta' \in \Theta_{ij}} I_{ij}(\theta, \theta') > 0$ for all $1 \leq i < I, 1 \leq j \leq J_i$ and $\theta \in \cup_{\ell > i} \Theta_\ell$.
A3. For each $j = 1, \ldots, J_i, i = 1, \ldots, I$ and $\theta \in \Theta$, $\{X_{ijt}, t \geq 0\}$ is a Markov chain on a state space $D$ with $\sigma$-algebra $\mathcal{D}$, irreducible with respect to a maximal irreducible measure on $(D, \mathcal{D})$ and aperiodic. Furthermore, $X_{ijt}$ is *Harris recurrent* in the sense that there exists a set $G_{ij} \in \mathcal{D}$, $\alpha_{ij} > 0$ and probability measure $\varphi_{ij}$ on $G_{ij}$ such that $P_{ij}^\theta\{X_{ijt} \in G_{ij} \text{ i.o.} | X_{ij0} = x\} = 1$ for all $x \in D$ and

$$(2.5)\ P_{ij}^\theta\{X_{ij1} \in A | X_{ij0} = x\} \geq \alpha_{ij} \varphi_{ij}(A) \quad \text{for all } x \in G_{ij} \text{ and } A \in \mathcal{D}.$$

A4. There exist constants $0 < \bar{b} < 1$, $b > 0$ and drift functions $V_{ij} : D \to [1, \infty)$ such that for all $j = 1, \ldots, J_i$ and $i = 1, \ldots, I$,

$$(2.6) \quad \sup_{x \in D} |g(x)|/V_{ij}(x) < \infty,$$

and for all $x \in D$, $\theta \in \Theta$,

$$(2.7) \quad P_{ij}^\theta V_{ij}(x) \leq (1 - \bar{b}) V_{ij}(x) + b \mathbf{1}_{G_{ij}}(x),$$

where $G_{ij}$ satisfies (2.5) and $P^{\theta}_{ij}V_{ij}(x) = \int_D V_{ij}(y)P^{\theta}_{ij}(x,dy)$. Moreover, we require that

$$(2.8) \quad \int_D V_{ij}(x)\nu_{ij}(dx;\theta)Q(dx) < \infty \quad \text{and} \quad V^*_{ij} := \sup_{x \in G_{ij}} V_{ij}(x) < \infty.$$

Let $\ell_{ij}(x,y;\theta,\theta') = \log[p_{ij}(x,y;\theta)/p_{ij}(x,y;\theta')]$ be the log likelihood ratio between $P^{\theta}_{ij}$ and $P^{\theta'}_{ij}$ and $N_\delta(\theta) = \{\theta' : \|\theta - \theta'\| < \delta\}$ a ball of radius $\delta$ around $\theta$, where $\|\cdot\|$ denotes Euclidean norm.

A5. There exists $\delta > 0$ such that for all $\theta, \theta' \in \Theta$,

$$(2.9) \quad K_{\theta,\theta'} := \sup_{x \in D} \frac{E_\theta[\sup_{\widetilde{\theta} \in N_\delta(\theta')} \ell^2_{ij}(X_{ij0}, X_{ij1}; \theta, \widetilde{\theta})|X_{ij0} = x]}{V_{ij}(x)} < \infty$$

for all $j = 1, \ldots, J_i$, $i = 1, \ldots, I$. Moreover,

$$(2.10) \quad \sup_{\widetilde{\theta} \in N_{\delta'}(\theta')} |\ell_{ij}(x,y;\theta',\widetilde{\theta})| \to 0 \text{ as } \delta' \to 0$$

for all $x, y \in D$ and $\theta' \in \Theta$.

Assumption A1 is a mild regularity condition to exclude unrealistic models. A2 is a positive information criterion: the first inequality makes sure that information is available in the first arm group to estimate $\theta$; while the second inequality allows us to collect information in the $i$-th arm group for moving to the next group when $\theta \in \Theta_\ell$ for some $\ell > i$. Assumption A3 is a recurrence condition and A4 is a drift condition. These two conditions are used to guarantee the stability of the Markov chain so that the strong law of large numbers and Wald's equation hold. A5 is a finite second moment condition that allows us to bound the probability that the MLE of $\theta$ lies outside a small neighborhood of $\theta$. This bound is important for us to determine the level of unequal allocation of observations that can be permitted in the testing stage of our procedure. The proof of the asymptotic lower bound in Theorem 1 requires only A1-A3; while additional A4 and A5 are required for the construction of efficient strategies attaining the lower bound.

## 2.2. Bad sets

The bad set is a useful concept for understanding the learning required within the group containing optimal arms. It is associated with the asymptotic lower bound described in Section 3 and is used explicitly in constructing the asymptotically efficient strategy. For $\theta \in \Theta_\ell$, define $J(\theta) = \{j : \mu^*(\theta) = \mu_{\ell j}(\theta)\}$ as the set of optimal jobs in group $\ell$. Hence $\theta \in \Theta_{\ell j}$ if and only if $j \in J(\theta)$. We also define the bad set, the set of 'bad' parameter values associated with $\theta$, as all $\theta' \in \Theta_\ell$ which cannot be distinguished from $\theta$ by processing any of the optimal jobs $\ell j$. Specifically,

$$(2.11) \quad B_\ell(\theta) = \Big\{\theta' \in \Theta_\ell \setminus \big(\bigcup_{j \in J(\theta)} \Theta_{\ell j}\big) : I_{\ell j}(\theta, \theta') = 0 \text{ for all } j \in J(\theta)\Big\}.$$

The bad set $B_\ell(\theta)$ is the intersection of two parameter sets. One set consists of parameter values that have different optimal arms from those for $\theta$. The other set contains parameter values that cannot be distinguished from sampling the optimal arm for $\theta$. When a parameter value is in the intersection, sampling from arms that are non-optimal for $\theta$ is required.

We note that if $I_{\ell j}(\theta, \theta') = 0$, then the transition probabilities of $X_{\ell j t}$ are identical under both $\theta$ and $\theta'$. If $\theta' \in B_\ell(\theta)$, then by definition, $\theta' \notin \cup_{j \in J(\theta)} \Theta_{\ell j}$ and hence $J(\theta') \cap J(\theta) = \emptyset$. Let $j \in J(\theta)$ and $j' \in J(\theta')$. Then $\mu_{\ell j'}(\theta') > \mu_{\ell j}(\theta') = \mu_{\ell j}(\theta) > \mu_{\ell j'}(\theta)$. Thus

(2.12) $$I_{\ell j'}(\theta, \theta') > 0 \text{ for all } \theta' \in B_\ell(\theta) \text{ and } j' \in J(\theta').$$

The interpretation of (2.12) is as follows. Although we cannot distinguish $\theta$ from $\theta' \in B_\ell(\theta)$ when sampling the optimal arm for $\theta$, we can distinguish them by sampling the optimal job for $\theta'$. This fact explains the necessity of processing non-optimal arms to collect information.

## 3. The regret lower bound

The following theorem gives an asymptotic lower bound for the regret (1.4) of uniformly good adaptive strategies under the partial order constraint $\preceq$. The proof can be found in [1]. We will discuss the relation of the lower bound with those in [6, 7] and [3].

**Theorem 1.** *Assume* A1-A3 *and let* $\theta \in \Theta_\ell$. *For any uniformly good adaptive strategy $\phi$ under the partial order constraint $\preceq$,*

(3.1) $$\liminf_{N \to \infty} R_N(\theta)/\log N \geq z(\theta, \ell),$$

*where $z(\theta, \ell)$ is the minimum value of the following minimization problem.*

(3.2) $$\text{Minimize} \sum_{i < \ell} \sum_{j=1}^{J_i} [\mu^*(\theta) - \mu_{ij}(\theta)] z_{ij}(\theta) + \sum_{j \notin J(\theta)} [\mu^*(\theta) - \mu_{\ell j}(\theta)] z_{\ell j}(\theta),$$
*subject to* $z_{ij}(\theta) \geq 0, \quad j = 1, \ldots, J_i, \text{ if } i < \ell, \ j \notin J(\theta), \text{ if } i = \ell,$

and

(3.3) $$\begin{cases} \inf_{\theta' \in \Theta_1}\{\sum_{j=1}^{J_1} I_{1j}(\theta, \theta') z_{1j}(\theta)\} \geq 1, \\ \inf_{\theta' \in \Theta_2}\{\sum_{j=1}^{J_1} I_{1j}(\theta, \theta') z_{1j}(\theta) + \sum_{j=1}^{J_2} I_{2j}(\theta, \theta') z_{2j}(\theta)\} \geq 1, \\ \vdots \\ \inf_{\theta' \in \Theta_{\ell-1}}\{\sum_{j=1}^{J_1} I_{1j}(\theta, \theta') z_{1j}(\theta) + \cdots + \sum_{j=1}^{J_{\ell-1}} I_{(\ell-1)j}(\theta, \theta') z_{(\ell-1)j}(\theta)\} \geq 1, \\ \inf_{\theta' \in B_\ell(\theta)}\{\sum_{i<\ell} \sum_{j=1}^{J_i} I_{ij}(\theta, \theta') z_{ij}(\theta) + \sum_{j \notin J(\theta)} I_{\ell j}(\theta, \theta') z_{\ell j}(\theta)\} \geq 1. \end{cases}$$

**Corollary 1.** *When there is only one group of arms, (3.1) reduces to the lower bound (1.11) of Lai and Robbins* [7].

*Proof.* When there is only group of arms, only the last inequality of (3.3) is needed and it takes the form

(3.4) $$\inf_{\theta' \in B(\theta)} \sum_{j \notin J(\theta)} I_j(\theta, \theta') z_j(\theta) \geq 1.$$

In [7], it is proved that

$$(3.5) \qquad E_\theta T_N(j) \geq \frac{\log N}{I(\theta_j, \theta^*)} \quad \text{for all } j \notin J(\theta),$$

where $\theta^* = \max_{1 \leq i \leq k} \theta_i$. Note that in [7], all jobs belong to the same family of probability distributions with different parameter values, and thus the KL information number does not depend on the job label but only the parameter value. Let $E_\theta T_N(j)/\log N = z_j(\theta)$, then (3.5) is the same as

$$(3.6) \qquad z_j(\theta) I(\theta_j, \theta^*) \geq 1 \quad \text{for all } j \notin J(\theta).$$

We first show that (3.4) $\Rightarrow$ (3.5). Because (3.4) implies that for all $\theta' \in B(\theta)$

$$(3.7) \qquad \sum_{j \notin J(\theta)} I(\theta_j, \theta'_j) z_j(\theta) \geq 1.$$

If $\theta' = (\theta'_1, \ldots, \theta'_k) \in B(\theta)$, then $\theta^* = \theta_{j^*} = \theta'_{j^*}$ and $\max_{1 \leq i \leq k} \theta'_i > \theta^*$. Suppose we choose a sequence of $\theta' \in B(\theta)$ such that there is only one component $\theta'_j$ approaching $\theta^*$ from above and other components $\theta'_{j'}$, $j' \notin J(\theta)$, all have the same values as the corresponding components of $\theta$. Taking infimum over this sequence of $\theta' \in B(\theta)$ in (3.7), we obtain (3.6). This complete the proof of (3.4) $\Rightarrow$ (3.5).

To prove (3.5) $\Rightarrow$ (3.4), we assume that (3.4) does not hold. That is, there exist a $\theta' \in B(\theta)$ such that

$$\sum_{j \notin J(\theta)} I(\theta_j, \theta'_j) z_j(\theta) < 1.$$

Because $\theta' \in B(\theta)$, there exists at least one component $\theta'_{j*}$ of $\theta'$ such that $\theta'_{j*} > \theta^*$. Then the preceding inequality and the property of exponential families imply that

$$z_{j^*}(\theta) I(\theta_{j^*}, \theta^*) < z_{j^*}(\theta) I(\theta_{j^*}, \theta'_{j^*}) < 1,$$

and thus (3.6) does not hold. This establishes (3.5) $\Rightarrow$ (3.4) and the proof is complete. $\square$

**Corollary 2.** *When there is only one arm in each group, then (3.1) reduces to the lower bound (1.17) of Hu and Wei [6].*

*Proof.* In Hu and Wei [6], the set $\Theta_i$ are intervals of $\Re$. Thus the infimum over $\Theta_i$ is achieved at the end points of the intervals. Furthermore, because there is only one arm in each group, the bad sets are all empty and therefore the last inequality in (3.3) is not needed. In view of these facts, it is straightforward to show that the systems of inequalities (3.3) reduces to (1.14) of Hu and Wei [6]. The proof is complete. $\square$

**Corollary 3.** *When there is only one arm in each group, the lower bound (3.1) reduces to (3.2) of Fuh and Hu [3].*

*Proof.* The assumptions A3 and A4 of Fuh and Hu [3] correspond to the regularity condition A1 and the positive information criterion A2 in Section 2, respectively. The A1, A2 and A5 of Fuh and Hu are essentially the same as Harris recurrence condition A3, the drift condition A4, and the finite second moment condition A5 of this paper, respectively.

Note that the definition of bad sets in [3] is different from that of this paper. In [3], the bad set consists of all those parameter values having optimal arm *not* in the same group and cannot be distinguished when sampling from the optimal arm. Here the bad set consists of parameter values that has different optimal arm (*but still in the same group*), and cannot be distinguished when sampling from the optimal arm(s). If we adopt the definition (2.11), then it is clear that the bad sets are all empty under the setting of [3].

The infimums in Problem A of Fuh and Hu [3] is taken over the union of $\Theta_i$ and the corresponding bad set. Because the bad sets in [3] are all empty as we point out earlier, the infimums is actually taken over $\Theta_i$. With this understanding, it is straightforward to verify that the lower bound (3.1) reduces to (3.2) of Fuh and Hu [3]. □

## 4. Super efficiency

The strategy in the allocation of the observations is as follows. For the rationale of the proposed strategy and more detailed discussion, please see [1]. Let $n_0$ and $n_1$ be positive integers that increase to infinity with respect to $N$ and satisfies $n_0 = o(\log N)$ and $n_1 = o(n_0)$.

1. *Estimation*. Select $n_0$ observations from each arm in group 1 and let $\widehat{\theta}$ be the maximum likelihood estimate (MLE) of $\theta$ defined by

$$(4.1) \qquad L(\theta) = \sum_{j=1}^{J} \sum_{t=1}^{n_0} \log p_{1j}(X_{1j(t-1)}, X_{1jt}; \theta), \quad \widehat{\theta} = \arg\max_{\theta \in \Theta} L(\theta).$$

Let $\ell = \min\{i : N_{\delta/2}(\widehat{\theta}) \cap \Theta_i \neq \emptyset\}$. Select an adjusted MLE estimate $\widehat{\theta}_a \in N_{\delta/2}(\widehat{\theta}) \cap \Theta_\ell$, (where $\delta \to 0$ as $N \to \infty$ at a rate to be specified in Theorem 1 below), in the following manner. Let $|\cdot|$ denote the number of elements in a finite set and

$$(4.2) \qquad \mathbf{J} = \max\{|J(\theta')| : \theta' \in N_{\delta/2}(\widehat{\theta}) \cap \Theta_\ell\}.$$

We require that

$$(4.3) \qquad \widehat{\theta}_a \in H := \{\theta \in N_{\delta/2}(\widehat{\theta}) \cap \Theta_\ell : |J(\theta)| = \mathbf{J}\}.$$

The motivation behind considering an adjusted MLE is to estimate $J(\theta)$ and the set $\Theta_i$ that $\theta$ belongs to consistently. This has implications in the experimentation phase. We note that if $|J(\theta)| > 1$, then $J(\widehat{\theta})$ need not be consistent for $J(\theta)$ and if $\Theta_i$ lies on $\Theta_i \setminus \Theta_i^*$ [see (2.2) and (2.4)], then $\widehat{\theta}$ need not be consistently inside $\Theta_i$. Conversely, the probability that $J(\widehat{\theta}_a) = J(\theta)$ and $\widehat{\theta}_a$ lying inside $\Theta_i$ tends to 1 as $N \to \infty$.

Let

$$(4.4) \qquad B_\ell(\theta; \delta) = \cup_{\theta' \in H} B_\ell(\theta')$$

and let $\{\widehat{z}_{ij}\}_{1 \leq i \leq \ell, 1 \leq j \leq J_i}$ minimize

$$(4.5) \qquad \sum_{i<\ell} \sum_{j=1}^{J_i} [\mu^*(\widehat{\theta}_a) - \mu_{ij}(\widehat{\theta}_a)] z_{ij} + \sum_{j \notin J(\widehat{\theta}_a)} [\mu^*(\widehat{\theta}_a) - \mu_{\ell j}(\widehat{\theta}_a)] z_{\ell j}$$

subject to the constraints

$$
(4.6) \quad \begin{cases} \inf_{\theta' \in \Theta_1} \{\sum_{j=1}^{J_1} I_{1j}(\widehat{\theta}_a, \theta')z_{1j}\} \geq 1, \\ \inf_{\theta' \in \Theta_2} \{\sum_{j=1}^{J_1} I_{1j}(\widehat{\theta}_a, \theta')z_{1j} + \sum_{j=1}^{J_2} I_{2j}(\widehat{\theta}_a, \theta')z_{2j}\} \geq 1, \\ \vdots \\ \inf_{\theta' \in \Theta_{\ell-1}} \{\sum_{j=1}^{J_1} I_{1j}(\widehat{\theta}_a, \theta')z_{1j} + \cdots + \sum_{j=1}^{J_{\ell-1}} I_{(\ell-1)j}(\widehat{\theta}_a, \theta')z_{(\ell-1)j}\} \geq 1, \\ \inf_{\theta' \in B_\ell(\widehat{\theta};\delta)} \{\sum_{i<\ell} \sum_{j=1}^{J_i} I_{ij}(\widehat{\theta}_a, \theta')z_{ij} + \sum_{j \notin J(\widehat{\theta}_a)} I_{\ell j}(\widehat{\theta}_a, \theta')z_{\ell j}\} \geq 1. \end{cases}
$$

Let $k = 1$.

2. *Experimentation.* If $k \leq \ell$, select $\lfloor \widehat{z}_{kj} \log N \rfloor$ observations from arm $kj$, where $\lfloor \cdot \rfloor$ denotes the greatest integer function. If $k > \ell$, we skip the experimentation stage. We note that if $B_\ell(\widehat{\theta}; \delta)$ is empty, then the last inequality in (4.6) is automatically satisfied and hence we can select $\widehat{z}_{\ell 1} = \cdots = \widehat{z}_{\ell J_\ell} = 0$. In other words, if $B_\ell(\widehat{\theta}; \delta)$ is empty, then the experimentation stage is also skipped over for $k = \ell$.

3. *Testing.* Start with a full set $\{k1, \ldots, kJ_k\}$ of unrejected jobs. The rejection of a job is based on the following test statistic. Let $F_k$, $1 \leq k \leq I$, be a probability distribution with positive probability on all open subsets of $\cup_{i=k}^{I} \Theta_i$. Define
(4.7)
$$
U_k(\mathbf{n}; \lambda) = \frac{\int_{\cup_{i=k}^{I} \Theta_i} \prod_{i=1}^{k} \prod_{j=1}^{J_i} \nu_{ij}(X_{ij0}; \theta) \prod_{t=1}^{n_{ij}} p_{ij}(X_{ij(t-1)}, X_{ijt}; \theta) \, dF_k(\theta)}{\prod_{i=1}^{k} \prod_{j=1}^{J_i} \nu_{ij}(X_{ij0}; \lambda) \prod_{t=1}^{n_{ij}} p_{ij}(X_{ij(t-1)}, X_{ijt}; \lambda)}
$$

for all $\lambda \in \Theta_k$.

(a) If $\widehat{\theta} \in \cup_{i>k} \Theta_i$: Add one observation from each unrejected job. Reject parameter $\lambda$ if $U_k(\mathbf{n}; \lambda) \geq N$. Reject a job $kj$ if all $\lambda \in \Theta_{kj}$ have been rejected at some point in the testing stage. If there is a job in group $k$ left unrejected and the total number of observations is less than $N$, repeat 3(a). Otherwise go to step 4.

(b) If $\widehat{\theta} \in \Theta_k$: Add $n_1$ observations from each unrejected job $kj$, $j \in J(\widehat{\theta})$ and one observation from each unrejected job $kj$, $j \notin J(\widehat{\theta})$. Reject a job $kj$ if all $\lambda \in \Theta_{kj}$ have been rejected at some point in the testing phase. If there is a job in group $k$ left unrejected and the total number of observations is less than $N$, repeat 3(b). Otherwise, go to step 4.

(c) If $\widehat{\theta} \in \cup_{i<k} \Theta_i$: Adopt the procedure of 3(a).

4. *Moving to the next group and termination.* The strategy terminates once $N$ observations have been collected. Otherwise, if $k < I$, increment $k$ by 1 and go to step 2; if $k = I$, select all remaining observations from a job $Ij$ satisfying $\mu_{Ij}(\widehat{\theta}) = \max_{1 \leq h \leq J_I} \mu_{Ih}(\widehat{\theta})$.

In [1] Theorem 2, it was established that when $B_\ell(\theta)$ is non-empty, then the asymptotic lower bound of the regret is attained with the procedure above. We shall show that the same procedure is not only asymptotically optimal but also the regret from the optimal group will be $o(\log N)$ when $B_\ell(\theta) = \emptyset$ as oppose to $O(\log N)$ when $B_\ell(\theta) \neq \emptyset$. An important key step required in our proof is the consistency result

(4.8) $$P_\theta\{B_\ell(\widehat{\theta}, \delta) = \emptyset\} \to 1 \text{ as } N \to \infty$$

under the empty bad set assumption.

**Theorem 2.** *Let $\theta \in \Theta_\ell$. Assume* A1-A5 *and* (1.5) . *Let $n_0 \to \infty$ with $n_0 = o(\log N)$ and $n_1 \to \infty$ such that $n_1 = o(n_0)$. There exists $\delta(= \delta_N) \downarrow 0$ as $N \to \infty$ such that*

$$(4.9) \qquad P_\theta\{\widehat{\theta} \in \Theta \setminus N_\delta(\theta)\} = o(n_1^{-1}) \text{ as } n_1 \to \infty.$$

*Moreover, if $B_\ell(\theta) = \emptyset$, then* (4.8) *holds and*

$$(4.10) \qquad \sum_{j=1}^{J_\ell} E_\theta T_N(\ell j) = o(\log N).$$

*Hence*

$$(4.11) \qquad \lim_{N \to \infty} R_N(\theta)/\log N = z(\theta, \ell).$$

*Proof.* The consistency of $\widehat{\theta}$ in (4.9) follows from A2 and (4.5) of Chan, Fuh and Hu [1]. We shall now prove (4.8). Since $\delta \downarrow 0$ and $\widehat{\theta}$ is consistent for $\theta$, it suffices from the definition of $B_\ell(\theta; \delta)$ in (4.4) to show that there exists $\delta_0 > 0$ such that

$$(4.12) \qquad B_\ell(\widetilde{\theta}) = \emptyset \text{ for all } \widetilde{\theta} \in N_{\delta_0}(\theta) \cap \Theta_\ell \text{ with } |J(\widetilde{\theta})| = \mathbf{J}.$$

We observe from the continuity of $\mu_{\ell j}$ that there exists $\delta_1 > 0$ such that $J(\widetilde{\theta}) \subset J(\theta)$ for all $\widetilde{\theta} \in N_{\delta_1}(\theta) \cap \Theta_\ell$. Hence it follows that if $|J(\widetilde{\theta})| = |J(\theta)|$, then it must be true that $J(\widetilde{\theta}) = J(\theta)$. We see from the definition of bad sets in (2.11) that for each $\theta' \in \Theta_\ell \setminus (\cup_{j \in J(\theta)} \Theta_{\ell j})$, $I_{\ell j}(\theta, \theta') > 0$ for some $j \in J(\theta)$ and hence by the continuity of the Kullback-Leibler information, there exists $\delta_2 > 0$ such that $I_{\ell j}(\widetilde{\theta}, \theta') > 0$ whenever $\widetilde{\theta} \in N_{\delta_2}(\theta)$. Select $\delta_0 = \min\{\delta_1, \delta_2\}$. Then (4.12) holds.

We shall next show (4.10). By (4.8) and since the experimentation stage is skipped over when $k = \ell$ and $B_\ell(\widehat{\theta}, \delta) = \emptyset$, it suffices to show that the expected total number of observations taken from inferior arms in the testing stage is $o(\log N)$. Define $p_N = P_\theta\{J(\widehat{\theta}_a) = J(\theta)\}$. Then by (4.3), (4.9) and as $J(\widetilde{\theta}) \subset J(\theta)$ for all $\widetilde{\theta} \in N_{\delta_1}(\theta)$ for some $\delta_1 > 0$, $1 - p_N = o(n_1^{-1})$. By (2.16) and the assumption $B_\ell(\theta) = \emptyset$, at least one optimal arm will provide positive information against each $\theta' \notin \cup_{j \in J(\theta)} \Theta_j$. By A3-A5 and (6.4), (6.5) of Chan, Fuh and Hu [1], (an expected) $O(\log N)$ number of observations from arms with positive information is required to reject each $\theta' \in \Theta_\ell \setminus (\cup_{j \in J(\theta)} \Theta_{\ell j})$. Hence $O(n_1^{-1} \log N)$ number of recursions is involved when $J(\theta) = J(\widehat{\theta}_a)$ because at least $n_1$ observations in each recursion has positive information. Similarly, $O(\log N)$ recursions is needed when $J(\theta) \neq J(\widehat{\theta}_a)$ because at least one observation in each recursion has positive information. The number of observations from inferior arms in each recursion is $O(1)$ if $J(\widehat{\theta}_a) = J(\theta)$ and $O(n_1)$ otherwise. Hence the expected number of observations from inferior arms during the recursion steps in the testing phase is

$$(4.13) \qquad p_N O(n_1^{-1} \log N) + (1 - p_N) O(n_1 \log N) = o(\log N).$$

The asymptotic result (4.11) follows from (4.10) and the proof of Chan, Fuh and Hu [1] Theorem 2. □

For the special case $\ell = 1$, it follows from (4.11) that $R_N(\theta) = o(\log N)$ occurs. In [2] and [10], a uniformly good procedure was proposed that satisfies $R_N(\theta) = O(1)$ when $\Theta$ is finite and $I = 1$.

## 5. The switching cost

Let $a(\theta) > 0$ be the switching cost between two arms and are not both optimal when the underlying parameter is $\theta$. It is assumed here that there is no switching cost when both arms are optimal. Then

$$L_N(\theta) := a(\theta) E_\theta \Big( \sum_{t=1}^{N-1} \mathbf{1}_{\{\phi_t \neq \phi_{t+1}, \min[\mu_{\phi_t}(\theta), \mu_{\phi_{t+1}}(\theta)] < \mu^*(\theta)\}} \Big)$$

is the average switching cost of a procedure. It is also desirable that this cost is asymptotically negligible compared to the regret as $N \to \infty$.

**Theorem 3.** *Under Assumptions A1 - A5, the strategy $\phi^*$ has average switching cost*

(5.1) $$L_N(\theta) = o(\log N) \text{ as } N \to \infty.$$

*Hence, the strategy is asymptotically optimal when there is switching cost.*

*Proof.* In the estimation stage it is require to take $n_0$ observations from each arm in group 1. We can take the $n_0$ observations in batches and switch only $J_1 - 1$ times. Therefore the switching cost from estimation stage is $a(\theta)(J_1 - 1)$. In the experimentation stage, we need to allocate at most $\widehat{z}_{kj} \log N$ observations to arm $kj$. Again this can be done in batches and thus the switching cost from experimentation stage is at most $a(\theta)(J_k - 1)$. In the testing stage, it is shown in (6.12) of Chan, Fuh and Hu [1], that the expected total number of observations is $o(\log N)$ and thus the switching cost is no more than $o(\log N)$. Adding the switching costs from the estimation, experimentation, and testing stages together, shows that the total cost due to switching is $o(\log N)$. However, the regret lower bound is $O(\log N)$, which implies that the switching cost constitutes a negligible part of the total regret as $n \to \infty$. This completes the proof that the proposed strategy is still asymptotically optimal with constant cost per switch. □

## 6. Extension of Wald's equation to Markovian rewards

As we will be focusing on a single arm $ij$ and fixed parameters $\theta_0$, $\theta_q$ such that $\mu := I_{ij}(\theta_0, \theta_q) > 0$ we will drop some of the references to $i$, $j$, $\theta_0$, $\theta_q$ and $q$ in this section. This applies also to the notations in assumptions A3-A5. Moreover, we shall use the notation $E(\cdot)$ as a short form of $E_{\theta_0}(\cdot)$ and $E_x(\cdot)$ as a short form of $E_{\theta_0}(\cdot | X_0 = x)$. Let $S_n = \xi_1 + \cdots + \xi_n$, where $\xi_k = \log[p_{ij}(X_{k-1}, X_k; \theta_0)/p_{ij}(X_{k-1}, X_k; \theta_q)]$ has stationary mean $\mu$ under $P_{\theta_0}$ and let $\tau$ be a stopping-time. We shall show that

(6.1) $$ES_\tau = \mu(E\tau) - E[\gamma(X_\tau)] + E[\gamma(X_0)]$$

for some function $\gamma$ to be specified in Lemma 1. In Lemma 2, we show that the conditions on $V$ in A4-A5 lead to bounds on $\gamma(x)$ and by applying Lemma 3, we obtain

(6.2) $$E|\gamma(X_\tau)| + E|\gamma(X_0)| = o(E\tau).$$

Substituting (6.2) back into (6.1), Wald's equation

(6.3) $$ES_\tau = [\mu + o(1)]E\tau$$

is established for Markovian rewards. Under uniform recurrence condition, Fuh and Lai [4] established Wald's equation based on perturbation theory for the transition operator. The Wald's equation was proved under the assumption that the solution for the Poisson equation exists in [5] based on Poisson equation for the transition operator. In this section, we apply the idea of regeneration epoch to derive the Wald's equation for Markov random walks.

By (2.5), we can augment the Markov additive process and create a split chain containing an atom, so that increments in $S_n$ between visits to the atom are independent. More specifically, we construct stopping-times $0 < \kappa(1) < \kappa(2) < \cdots$ using an auxiliary randomization procedure such that

$$(6.4) \quad P\{X_{n+1} \in A, \kappa(i) = n+1 | X_n = x, \kappa(i) > n \geq \kappa(i-1)\} = \begin{cases} \alpha\varphi(A) & x \in G, \\ 0 & otherwise. \end{cases}$$

Then by Lemma 3.1 of Ney and Nummelin [9],

(i) $\{\kappa(i+1) - \kappa(i) : i = 1, 2, \ldots\}$ are i.i.d. random variables.
(ii) the random blocks $\{X_{\kappa(i)}, \ldots, X_{\kappa(i+1)-1}\}$, $i = 1, 2, \ldots$, are independent and
(iii) $P\{X_{\kappa(i)} \in A | \mathcal{F}_{\kappa(i)-1}\} = \varphi(A)$, where $\mathcal{F}_n = \sigma$-field generated by $\{X_0, \ldots, X_n\}$.

By (ii)-(iii), $E_\varphi(S_\kappa - \kappa\mu) = 0$. Define $\kappa = \kappa(1)$. We shall use the notation "$n =$ atom" to denote $n = \kappa(i)$ for some $i$.

**Lemma 1.** *Let $\gamma(x) = E_x(S_\kappa - \kappa\mu)$. Then $Z_n = (S_n - n\mu) + \gamma(X_n)$ is a martingale with respect to $\mathcal{F}_n$. Hence (6.1) holds.*

*Proof.* We can express

$$Z_n = E(S_{U_n} - U_n\mu | \mathcal{F}_n) \quad \text{where } U_n = \inf\{m > n : m = \text{atom}\}.$$

If $X_n = x_n \notin G$, then by (6.4), $U_n > n+1$. Hence $U_{n+1} = U_n$ and

$$(6.5) \quad E(Z_{n+1} | \mathcal{F}_n) = Z_n$$

because $\mathcal{F}_{n+1} \supset \mathcal{F}_n$. If $X_n = x_n \in G$, then by (6.4) and (ii),

$$E(Z_{n+1} | \mathcal{F}_n) - Z_n = E[(S_{U_{n+1}} - S_{U_n}) + (U_{n+1} - U_n) | \mathcal{F}_n] = \alpha E_\varphi(S_\kappa - \kappa\mu) = 0$$

and hence (6.5) also holds. □

**Lemma 2.** *Under* A3-A5,

$$|\gamma(x)| \leq \beta^{-1}[V(x) + b + (V^* + b)V^*(\alpha^{-1} + 1)](K + 1 + |\mu|),$$

*where $\alpha$ satisfies (2.5), $V^*$ is defined in* A4 *and $K$ is defined in (2.9).*

*Proof.* By (2.9),

$$(6.6) \quad V(x) \geq K^{-1}E_x\xi_1^2 \geq K^{-1}(E_x|\xi_1| - 1).$$

Let $0 < \sigma(1) < \sigma(2) < \cdots$ be the hitting times of the set $G$ and let $\sigma = \sigma(1)$. Let

$$(6.7) \quad m_n(A) = E_x\Big[\sum_{n=1}^{\kappa} V(X_n)\mathbf{1}_{\{X_n \in A\}}\Big]$$

for all measurable set $A \subset D$. By (2.7),

$$E_x[V(X_n)\mathbf{1}_{\{\sigma \geq n\}}] \leq (1-\beta)E_x[V(X_{n-1})\mathbf{1}_{\{\sigma \geq n-1\}}], \ n \geq 2$$

and

$$E_x[V(X_1)] \leq V(x) + b.$$

Hence by induction,

$$(6.8) \qquad E_x\Big[\sum_{n=1}^{\sigma} V(X_n)\Big] \leq [V(x) + b]\sum_{n=1}^{\infty}(1-\beta)^{n-1} = [V(x) + b]/\beta.$$

By (6.7)-(6.8), and as $V \geq 1$,

$$m_n(D) = E_x\Big\{\sum_{n=1}^{\sigma} V(X_n) + \sum_{k=1}^{\infty}\Big[\sum_{n=\sigma(k)+1}^{\sigma(k+1)} V(X_n)\Big]\mathbf{1}_{\{\kappa > \sigma(k)\}}\Big\}$$

$$= E_x\Big\{\sum_{n=1}^{\sigma} V(X_n) + \sum_{k=1}^{\infty} E_{X_{\sigma(k)}}\Big[\sum_{n=1}^{\sigma} V(X_n)\Big]\mathbf{1}_{\{\kappa > \sigma(k)\}}\Big\}$$

$$\leq \beta^{-1}[V(x) + b] + E_x\Big\{\sum_{k=1}^{\infty} \beta^{-1}[V(X_{\sigma(k)}) + b]\mathbf{1}_{\{\kappa > \sigma(k)\}}\Big\}$$

$$(6.9) \qquad \leq \beta^{-1}[V(x) + b] + \beta^{-1}(V^* + b)m_n(G).$$

But by (6.4), $m_n(G) \leq V^*(\alpha^{-1} + 1)$. Since $\gamma(x) \leq (K + 1 + |\mu|)m_n(D)$, Lemma 2 holds. □

Let $W_i = |\gamma(X_{\kappa(i)})| + \cdots + |\gamma(X_{\kappa(i+1)-1})|$, for $i \geq 1$. Then by A3-A5, Lemma 4 and its proof, and (i)-(iii), $W_1, W_2, \ldots$ are i.i.d. with finite mean while by (2.8), $W_0 := |\gamma(X_0)| + \cdots + |\gamma(X_{\kappa(1)-1})|$ also has finite mean.

**Lemma 3.** *Let $M_n = \max_{1 \leq k \leq n} W_k$. Then for any stopping-time $\tau$,*

$$(6.10) \qquad E(M_\tau) = o(E\tau).$$

*Proof.* Let $\delta > 0$ and let $c(= c_\delta) > 0$ be large enough such that $E[(W_1 - c)^+] \leq \delta$. We shall show that

$$Z_n = (M_n \vee c) - n\delta$$

is a supermartingale. Indeed for any $\lambda \geq 0$,

$$E[M_{n+1} \vee c | M_n \vee c = c + \lambda] = c + \lambda + E[(W_{n+1} - c - \lambda)^+] \leq c + \lambda + \delta$$

and the claim is shown. Hence $EZ_\tau \leq EZ_0 = c$ and it follows that $E(M_\tau) \leq E(M_\tau \vee c) \leq \delta(E\tau) + c$. Lemma 3 then follows by letting $\delta \downarrow 0$. □

## 7. Appendix

*Proof.* Proof of (1.3) Let $X_{ijt}$ denotes the $t$th observation taken from arm $ij$. Then

$$(7.1) \qquad \Big|W_N(\theta) - \sum_{i=1}^{I}\sum_{j=1}^{J_i} \mu_{ij}(\theta)E_\theta T_N(ij)\Big| \leq \sum_{i=1}^{I}\sum_{j=1}^{J_i}\sum_{t=1}^{\infty} |E_\theta g(X_{ijt}) - \mu_{ij}(\theta)|.$$

For any signed measure $\lambda$ on $(D, \mathcal{D})$, let

$$\|\lambda\|_{V_{ij}} = \sup_{h:|h|\leq V_{ij}} \left| \int h(x)\lambda(dx) \right|. \tag{7.2}$$

It follows from Meyn and Tweedie ([8], p.367 and Theorem 16.0.1) that under A3 and the geometric drift condition (2.7),

$$\omega_{ij} := \sup_{\theta \in \Theta, x \in D} \sum_{t=1}^{\infty} \|P_{ijt}^{\theta}(x, \cdot) - \pi_{ij}(\theta)\|_{V_{ij}}/V_{ij}(x) < \infty, \tag{7.3}$$

where $P_{ijt}^{\theta}(x, \cdot)$ denotes the distribution of $X_{ijt}$ conditioned on $X_{ij0} = x$ and $\pi_{ij}(\theta)$ denotes the stationary distribution of $X_{ijt}$ under parameter $\theta$. By (2.6), there exists $\kappa > 0$ such that $\kappa|g(x)| \leq V_{ij}(x)$ for all $x \in D$ and hence it follows from (7.2) and (7.3) that

$$\kappa \sum_{t=1}^{\infty} |E_{\theta,x} g(X_{ijt}) - \mu_{ij}(\theta)| \leq \omega_{ij} V_{ij}(x), \tag{7.4}$$

where $E_{\theta,x}$ denotes expectation with respect to $P_\theta$ and intial distribution $X_{ij0} = x$.

In general, for any initial distribution $\nu_{ij}(\cdot; \theta)$, it follows from (2.8) and (7.4) that

$$\sum_{t=1}^{\infty} |E_\theta g(X_{ijt}) - \mu_{ij}(\theta)| \leq \int \sum_{t=1}^{\infty} |E_{\theta,x} g(X_{ijt}) - \mu_{ij}(\theta)| \nu_{ij}(x; \theta) Q(dx) < \infty$$

uniformly over $\theta \in \Theta$ and hence (1.3) follows from (7.1). □